\mag=\magstep1          
\documentstyle{amsppt}\input amsppt1
\pageheight{24true cm}\pagewidth{16true cm}
\parindent=4mm\parskip=3pt plus1pt minus.5pt
\nologo   
\NoBlackBoxes 
\NoRunningHeads

\def\i{\looparrowright}   
\def\e{\hookrightarrow}

\def\f{\flushpar }
\def\nl{\newline }

\def\x{\times }

\topmatter
\title
Link cobordism and the intersection of slice discs
\endtitle
\author
Eiji Ogasa\\
This paper is published in\\
Bulletin of the London Mathematical Society, 31, 1999, 1-8.\\
This manuscript is not the published version. 
\endauthor
\thanks{
This research was partially supported by Research Fellowships
of the Promotion of Science for Young Scientists.}
\endtopmatter
\document
\baselineskip11pt

Abstract. It is well-known that all 2-knots are slice. 
Are all 2-links slice? This is an outstanding open question. 
In this paper we prove the following:  
For any 2-component 2-link $(J,K)\subset S^4=\partial B^5\subset B^5$, 
there is an embedded disc $D^2_J$ (respectively, $D^2_K$) $\subset B^5$ 
with the following properties:  
$J$ (respectively $K$) bounds $D^2_J$ (respectively, $D^2_K$). 
 $D^2_J$ and $D^2_K$ intersect transversely. 
 $D^2_J\cap D^2_K$ in $D^2_J$ (respectively, $D^2_K$) 
is a trivial 1-knot.

\head 1. Definitions\endhead

We work in the smooth category. 
     
An {\it (oriented) (ordered) $m$-component n-(dimensional) link}
 is a smooth, oriented submanifold $L=\{K_1,\dots ,K_m\}$ of $S^{n+2}$, 
  which is the ordered disjoint union of $m$ manifolds, each PL homeomorphic 
to the standard $n$-sphere. If $m=1$, then $L$ is called a {\it knot}.

We say that $m$-component $n$-dimensional links, 
$L_0$ and $L_1$,  are {\it (link-)concordant} or {\it  (link-)cobordant} 
if there is a smooth oriented submanifold  
$\widetilde{C}$=\{$C_1$ ,\dots ,$C_m$\} of $S^{n+2}\times [0,1]$,  
which meets the boundary transversely in $\partial \widetilde{C}$, 
 is PL homeomorphic to  $L_0 \times[0,1]$ 
  and meets $S^{n+2}\times \{l\}$ in $L_l$ ($l=0,1$).
If $m=1$, we say that $n$-knots $L_0$ and $L_1$ are 
{\it (knot-)concordant} or {\it  (knot-)cobordant}.
Then we call $C$ a {\it concordance-cylinder}
 of  the two $n$-knots $L$ and $L'$.

If an $n$-link $L$ is concordant to the trivial link, then we call $L$ a {\it slice} link. 

If an $n$-link $L=\{K_1,\dots ,K_m\}$ $\subset S^{n+2}=\partial B^{n+3}\subset B^{n+3}$ is slice, 
then 
there is a disjoint union of $(n+1)$-discs, 
$D^{n+1}_1\amalg\dots\amalg D^{n+1}_m$  
in $B^{n+3}$ such that 
$D^{n+1}_i\cap S^{n+2}$= $\partial D^{n+1}_i$=$K_i$.     
$(D^{n+1}_1,\dots, D^{n+1}_m)$ is called  
{\it a set of  slice discs  for $L$}. 
If $m=1$, then  $D^{n+1}_1$ is called a {\it slice disc} for the knot $L$.

Note. Let $L=(K_1,K_2)$ be an $n$-link. 
If we say $(D_1, D_2)$ is a set of slice disc for $L$, 
then we assume $D_1\cap D_2=\phi$. 
If we say that $D_1$ is a slice disc for $K_1$ 
and that $D_2$ is a slice disc for $K_2$, 
then  we do not assume $D_1\cap D_2=\phi$.

\head 2. Open problems \endhead

We have the following outstanding open problems.
See [1] for results and history arround them. 

{\bf Problem 2.1}  
Let $L_1$ and $L_2$ be $m$-component $n$-links ($n\geqq1$).
Then are $L_1$ and $L_2$ link concordant?

In particular, we have the following problems. 

{\bf Problem 2.2}
Is every even dimensional link slice?

{\bf Problem 2.2.1}
Is every 2-dimensional link slice?

{\bf Problem 2.3}
Is every odd dimensional link concordant to (a sublink of) a 
homology boundary link?   
( See [1] 
for homology boundary links.)

\head 3. Our problems \endhead

We modify Problem 2.2.1 to formulate the following Problem 3.1. 

We consider the case of 2-component links.
Let $L=(K_1, K_2)$ $\subset S^4$ $=\partial B^5$ $\subset B^5$  
be a $2$-link.  
By Kervaire's theorem in [2] 
there exist $3$-discs $D^3_i$ ($i=1,2$) embedded in $B^5$  
such that $D^3_i$ $\cap$ $S^4$=$\partial D^3_i$=$K_i$. 
Then  $D^3_1$ and $D^3_2$ intersect mutually in general. 
Furthermore, $D^3_1$ $\cap$ $D^3_2$ in  $D^3_i$ 
defines a $1$-link ($i=1,2$).

{\bf Problem 3.1}
Can we remove the above intersection $D^3_1$ $\cap$ $D^3_2$ 
by modifying an embedding of $D^3_1$ and  an embedding of $D^3_2$?

The answer to Problem 3.1 is affirmative if and only if   
the answer to Problem 2.2.1 is affirmative. 

Here we make another problem in order to attack Problem 3.1.   

{\bf Problem 3.2}
What is obtained as a pair of  1-links 
($D^3_1$ $\cap$ $D^3_2$ in  $D^3_1$,    
 $D^3_1$ $\cap$ $D^3_2$ in  $D^3_2$)   
by modifying an embedding of $D^3_1$ and  an embedding of $D^3_2$?

In this paper we give some answers to Problem 3.2 and its high dimensional version. 

\head 4. Main results \endhead

We have the following theorems, 
which are answers to Problem 3.2 and its high dimensional verion. 

\proclaim{Theorem 4.1}
For all 2-component $2m$-links  $L=(K_1, K_2)$ 
$\subset S^{2m+2}$ $=\partial B^{2m+3}$ $\subset B^{2m+3}$  
($2m\geq2$),  
there exist $(2m+1)$-disics $D^{2m+1}_1$ and  $D^{2m+1}_2$ $\subset B^{2m+3}$
 such that 
$D^{2m+1}_i$ is a slice disc for $K_i$ and 
$D^{2m+1}_1$ $\cap$ $D^{2m+1}_2$ in  $D^{2m+1}_i$ 
is the trivial $(2m-1)$-knot  ($i=1,2$).   
 \endproclaim

In particular, we have the following, which is the $2m=2$ case of Theorem 4.1.
\proclaim{Theorem 4.1.1}
Let   $L=(K_1, K_2)$ be a 2-component 2-link.   
Then there are 3-discs, $D^3_1$  and $D^3_2$, as in \S3 
such that  
each 1-link 
$D^3_1$ $\cap$ $D^3_2$ in  $D^3_i$ ($i=1,2$) is the trivial 1-knot. 
 \endproclaim

All even dimensional boundary links are slice. 
Hence the following Theorem 4.2 is Theorem 4.1 when $n$ is even.

\proclaim{Theorem 4.2}
For all 2-component $n$-link  $L=(K_1, K_2)$ ($n>1$),  
there exist  a boundary link  $L'=(K'_1, K'_2)$  satisfying that  
$K'_i$ is concordant to $K_i$ and a concordance-cylinder 
$\cases
\text{$C_1$  } \\ 
\text{$C_2$  } \\ 
\endcases$
of 
$\cases
\text{$K_1$  } \\ 
\text{$K_2$  } \\ 
\endcases$
and  
$\cases
\text{$K'_1$  } \\ 
\text{$K'_2$  } \\ 
\endcases$   
 such that 
$C_1$ $\cap$ $C_2$ in $C_i$ is the trivial $(n-1)$-knot ($i=1,2$).  
\endproclaim

By [1] 
not all 2-component odd dimensional links are concordant to boundary links. 
Therefore, when $n$ is odd, Theorem 4.2 is best possible from a viewpoint. 

Note. In this paper we do not solve Problem 2.1, 2.2, 2.3, 3.1 completely. 
They are now open problems.

\head 5. Proof of Theorem 4.1.1 \endhead

Let $L=(K_a, K_b)$ be a 2-link in $S^4$. 
There exists a Seifert hypersurface $V_i$ for $K_i$ $(i=a,b)$ 
such that $V_a\cap K_b$=$V_b\cap K_a$=$\phi$ 
and that $V_a$ intersects $V_b$ transversely. (See e.g. [7].)

Then $F$= $V_a\cap V_b$ is a closed oriented surface. 
We give $F$ a spin structure $\sigma$ naturally 
by using  $V_a$, $V_b$ and $S^4$.   

\proclaim{Claim} 
We can suppose that $F$ is connected. 
\endproclaim

Proof. 
Put $V_a\cap V_b=$ $F$=$F_1\amalg...\amalg F_\nu (\nu\geqq2)$.  
It suffices to prove that thare is a Seifert hypersurface $V'_b$ 
such that the connnected components of $V_a\cap V'_b$ are $\nu-1$.  

We made $V'_b$ from $V_b$ by using an embedded 1-handles as follows. 


Take points $p\in F_1$ and $q\in F_2$.
Let $l$ be a path in $V_a$ which connects $p$ and $q$. 
Take a 3-dimensional 1-handle $h$ $(\subset V_a)$ whose core is $l$ 
and which is attached to $F_1$ and $F_2$. 
Call the attach part, which is 2-discs, $P$ and $Q$. 
Recall $P,Q$ $\subset V_b$. 
Let $\widetilde h$ be a tubular neighborhood of $h$ in $S^4$.  
Let $\widetilde h$ be a 4-dimensional 1-handle which is attached to $V_b$. 
Let $\widetilde P$ (resp. $\widetilde Q$) be a tubular neighborhood of 
$P$ (resp. $Q$) in $V_b$, which is 3-discs. 
Then the attach part of  $\widetilde h$ is   $\widetilde P\amalg$  $\widetilde Q$. 
We carry out the surgery on  $V_b$ by using $\widetilde h$ and 
make a 3-manifold from $V_b$.  Call it $V'_b$.  
Then $V_a\cap V'_b$= $(F_1\sharp F_2)\amalg F_3\amalg...\amalg F_\nu$. 
Hence the connnected components of $V_a\cap V'_b$ are $\nu-1$.  

By the induction, Claim holds.

Consider the Sato-Levine invariant $\beta(L)$ of $L$. 
It is defined in [7]. 
$\beta(L)$=$[(F, \sigma)]$ $\in$ $\Omega^{\roman{spin}}_2\cong\Bbb Z_2$ 
by [5].
By [6],  
$\beta(L)=0$.

There exist simple closed curves 
$C_1,...,C_l$,$C'_1,...,C'_l$
 in $F$, 
where $l$ is the genus of $F$,   
such that the set of 1-cycles 
$[C_1],...,[C_l]$,$[C'_1],...,[C'_l]$  
represents a symplectic basis of $H_1(F:\Bbb Z)$.    
By P.36 of [3], we can suppose 
$[(C_i,\sigma\vert C_i)]$=0$\in$ $\Omega^{\roman{spin}}_1\cong\Bbb Z_2$
for each $i$.

\proclaim{Lemma 1}
Let  $S^2_a$ and $S^2_b$ be diffeomorphic to the 2-sphere. 

There exists an embedding 
$f:S^2_a\x[0,1]\amalg S^2_b\x[0,1]\e X$ with the following properties. 

\roster
\item
$X\cong\overline{(\sharp^p S^2\x B^3)-B^5}$. 
Put $\partial X$=$S^4_0\amalg M$. 
\nl
$S^4_0$ is diffeomorphic to the standard 4-sphere. 
$M$ is diffeomorphic to $\sharp^p S^2\x S^2$.
\item
$f$ is transverse to $S^4_0$ and 
$f(S^2_i\x[0,1])\cap S^4_0$=$f(S^2_i\x\{0\})$ ($i=a,b$).
\nl
$f(S^2_a\x\{0\}\amalg S^2_b\x\{0\})$ in $S^4_0$ is the 2-link $L$. 

\item
$f$ is transverse to $M$ and 
$f(S^2_i\x[0,1])\cap M$=$f(S^2_i\x\{1\})$ ($i=a,b$).
\nl
$f(S^2_i\x\{1\})$=$\partial\widetilde{V_i}$ 
$\subset$ $\widetilde{V_i}$ $\subset$ $M$  
and 
$\widetilde{V_i}$ is 
a compact connected oriented codimension one submanifold of $M$.
\nl
$\widetilde{V_a}$ and $\widetilde{V_b}$ intersect transversely   
and 
$\widetilde{V_a}$ $\cap$ $\widetilde{V_b}$ is diffeomorphic to the 2-sphere. 
\endroster
\endproclaim

\f{\bf Proof of Lemma 1}

Let 
$\cases
\text{$N_F(C_i)$} \\
\text{$N_a(C_i)$} \\
\text{$N_b(C_i)$} \\
\text{$N_S(C_i)$} \\
\endcases$
be the tubular neighborhood of $C_i$ in 
$\cases
\text{$F$} \\
\text{$V_a$} \\
\text{$V_b$} \\
\text{$S^4$} \\
\endcases$. 

Then we can consider 
$\cases
\text{$N_a(C_i)=N_F(C_i)\x I_a\x \{0\}$} \\
\text{$N_b(C_i)=N_F(C_i)\x \{0\}\x I_b$} \\
\text{$N_S(C_i)=N_F(C_i)\x I_a\x I_b$}    \\
\endcases$, 

where $I_a$ and $I_b$ is diffeomorphic to [-1,1]. 
 
Let 
$\cases
\text{$h^2_3$} \\
\text{$h^2_a$} \\
\text{$h^2_b$} \\
\text{$h^2_5$} \\
\endcases$   
be the 
$\cases
\text{3-dimensional} \\
\text{4-dimensional} \\
\text{4-dimensional} \\
\text{5-dimensional} \\
\endcases$  2-handles. 

Let $D^2$ be diffeomorphic to the 2-disc. 
Let $J$, $J_a$ and $J_b$ be diffeomorphic to [-1,1].

We can put 
$\cases
\text{$h^2_5$=$D^2\x J\x J_a\x J_b$} \\
\text{$h^2_a$=$D^2\x J\x J_a\x \{0\}$} \\
\text{$h^2_b$=$D^2\x J\x\{0\}\x J_b$} \\
\text{$h^2_3$=$D^2\x J\x\{0\}\x\{0\}$} \\
\endcases$.

Take $S^4\x[0,1]$. 
We can suppose naturally that 
$V_a\x[0,1]$, 
$V_b\x[0,1]$, and 
$F\x[0,1]$ is embedded in $S^4\x[0,1]$ 
and that 
the each embedding is level preserving. 

We attach $h^2_3$ to $F\x[0,1]$ along $N_F(C_i)$ 
($\subset F\x\{1\}$$\subset S^4\x\{1\}$). 

$[(C_i,\sigma\vert C_i)]$=0$\in$ $\Omega^{\roman{spin}}_1$.
Hence the spin structure on $F\x[0,1]$ extends to $h^2_3$. $\cdot\cdot\cdot\cdot(\sharp)$

Then we can attach 
$\cases
\text{$h^2_a$} \\
\text{$h^2_b$} \\
\text{$h^2_5$} \\
\endcases$   
to 
$\cases
\text{$V_a\x[0,1]$} \\
\text{$V_b\x[0,1]$} \\
\text{$S^4\x[0,1]$} \\
\endcases$   
along
$\cases
\text{$N_a(C_i)$ ($\subset V_a\x\{1\}$)} \\
\text{$N_b(C_i)$ ($\subset V_b\x\{1\}$)} \\
\text{$N_S(C_i)$ ($\subset S^4\x\{1\}$)} \\
\endcases$
naturally. 

By $(\sharp)$, 
the spin structure on 
$\cases
\text{$V_a\x[0,1]$} \\
\text{$V_b\x[0,1]$} \\
\text{$S^4\x[0,1]$} \\
\endcases$   
extends to 
$\cases
\text{$h^2_a$} \\
\text{$h^2_b$} \\
\text{$h^2_5$} \\
\endcases$   
$\cdot\cdot\cdot\cdot(*)$.

For each $C_i$ ($i=1,...,l$), 
we attach 2-handles to  
$\cases
\text{$F\x[0,1]$} \\
\text{$V_a\x[0,1]$} \\
\text{$V_b\x[0,1]$} \\
\text{$S^4\x[0,1]$} \\
\endcases$   
as above. 

Put $\widetilde X$=$S^4\x[0,1]$ $\cup$ (the 5-dimensional 2-handles).

By $(*)$, the attaching maps of the 5-dimensional handles are spin preserving diffeomorphism. 
Therefore 
 $\widetilde X$=$\overline{\sharp^p S^2\x B^3-B^5}$
  for a non-positive integer $p$.

Put $X$= $\widetilde X$ and take $f$ to be 
$f(S^2_i\x[0,1])$=$(\partial V_i)\x[0,1]$ ($\subset S^4\x[0,1]$) ($i=a,b$). 

The 4-manifold made from $V_i$ by the surgeries is called $\widetilde{V_i}$

We can make the corners of the new manifolds smooth, if necessary. 

Then Lemma 1 holds by the construction of them. 

This completes the proof of Lemma 1.

Under the condition of Lemma 1, 
we have the following claim. 

\proclaim {Claim 2}
In Lemma 1, we can move each of 
$\widetilde{V_a}$ and 
$\widetilde{V_b}$ 
by isotopy so that the following conditions  hold. 
(Note. 
The image of new map is not homeo to the old image $\widetilde{V_a}\cup\widetilde{V_b}$.)

\roster
\item
$\widetilde{V_a}$ $\cap$ $\widetilde{V_b}$ 
is diffeomorphic to the annulus. 
Put $\partial $($\widetilde{V_a}$ $\cap$ $\widetilde{V_b}$)
=$S^1_a\amalg S^1_b$. 
\item
$S^1_a\subset f(S^2_a\x\{1\})$. 
$S^1_b\subset f(S^2_b\x\{1\})$. 
\item
The deformation moves $\partial\widetilde{V_a}\amalg \partial\widetilde{V_b}$ 
by isotopy. 
\endroster
\endproclaim

{\bf Proof of Claim 2.} 
We move $\widetilde{V_a}$ and $\widetilde{V_b}$. 

In Lemma 1, 
$\widetilde{V_a}\cap$ $\widetilde{V_b}$ is a 2-sphere $S^2$.  
Take points $p_a$, $p_b$ $\in S^2$. Let $p_a\neq p_b$. 
Take a path  
$\cases
\text{$l_a$} \\
\text{$l_b$} \\
\endcases$   
which connects 
$\cases
\text{$p_a$} \\
\text{$p_b$} \\
\endcases$   
and 
a point in 
$\cases
\text{$\partial\widetilde{V_a}=f(S^2_a)$} \\
\text{$\partial\widetilde{V_b}=f(S^2_b)$} \\
\endcases$

Note 
$l_a \subset \widetilde{V_a}$ and 
$(l_a-p_a) \cap \widetilde{V_b}=\phi$.  
Note    
$l_b \subset \widetilde{V_b}$ and 
$(l_b-p_b) \cap \widetilde{V_a}=\phi$.  

Take a neighborhood of $l_a$ in $\widetilde{V_a}$. Call it $N_a$. 
Suppose $N_a$ is diffeomorphic to the close 2-disc and $l_a\subset$ Int $N_a$. 

Take a neighborhood of $l_b$ in $\widetilde{V_b}$. Call it $N_b$. 
Suppose $N_b$ is diffeomorphic to the close 2-disc and $l_b\subset$ Int $N_b$.

We can move $\widetilde{V_a}$ by isotopy and coincide with  
 $\overline{\widetilde{V_a}-N_a}$. 
Note. When we move $\widetilde{V_a}$ by isotopy, 
$\partial\widetilde{V_a}=f(S^2_a)$ intersect $\widetilde{V_b}$  
but does not intersect $\partial\widetilde{V_b}=f(S^2_b)$.    

We can move $\widetilde{V_b}$ by isotopy and coincide with  
 $\overline{\widetilde{V_b}-N_b}$. 
Note. When we move $\widetilde{V_b}$ by isotopy, 
$\partial\widetilde{V_b}=f(S^2_b)$ intersect $\widetilde{V_a}$ 
but does not intersect $\partial\widetilde{V_a}=f(S^2_a)$.    

This completes the proof of Claim 2.

In particular, the following holds in Claim 2. 

\proclaim{Claim 2.1}
$(\widetilde{V_a}$ $\cap$ $\widetilde{V_b})$ $-S^1_a$ is 
in the interior of $\widetilde{V_a}$.   

$(\widetilde{V_a}$ $\cap$ $\widetilde{V_b})$ $-S^1_b$ is 
in the interior of $\widetilde{V_b}$.   
\endproclaim

\proclaim{Lemma 3}
Let $S^2_a$ and $S^2_b$ be diffeomorphic to the 2-sphere. 

There exists an embedding 
$g:S^2_a\x[0,1]\amalg S^2_b\x[0,1]\e Y$ with the following properties. 

\roster
\item
$Y\cong\overline{(\natural^q S^2\x B^3)-B^5}$. 
Put $\partial Y$=$S^4_0\amalg N$. 
\nl
$S^4_0$ is diffeomorphic to the standard 4-sphere. 
$N$ is diffeomorphic to $\sharp^q S^2\x S^2$.
\item
$g$ is transverse to $S^4_0$, and 
$g(S^2_i\x[0,1])\cap S^4_0$=$g(S^2_i\x\{0\})$ ($i=a,b$).
\nl
$g(S^2_a\x\{0\}\amalg S^2_b\x\{0\})$ in $S^4_0$ is the 2-link $L$. 
\item
$g$ is transverse to $N$, and 
$g(S^2_i\x[0,1])\cap N$=$g(S^2_i\x\{1\})$ ($i=a,b$).
\nl
$g(S^2_i\x\{1\})$=$\partial D^3_i$ $\subset$ $D^3_i$ $\subset$ $N$. 
$D^3_i$ is diffeomorphic to the 3-disc.
\nl
$D^3_a$ and $D^3_b$ intersect transversely, and 
$D^3_a$ $\cap$ $D^3_b$ is diffeomorphic to  the annulus. 
\nl 
Put $\partial$($D^3_a$ $\cap$ $D^3_b$) =$S^1_a\amalg S^1_b$. 
\nl
$S^1_a\subset g(S^2_a\x\{1\})$.
$S^1_b\subset g(S^2_b\x\{1\})$. 
\endroster
\endproclaim 

\f{\bf Proof of Lemma 3}
Take an oriented 3-manifold which is orientation-preserving diffeomorphic to $\widetilde{V_i}$ 
($i=a,b$).
Call it  $\widetilde{V_i}$ again. 

Let $\widetilde{V_i}\cup D^3$ denote 
the oriented connected closed 3-manifold 
made from $\widetilde{V_i}$ and the 3-disc $D^3$. 
There exists a spin structure on  $\widetilde{V_i}$ induced from 
the unique spin structure on $S^4$.
The spin structure extends over  $\widetilde{V_i}\cup D^3$ uniquely.  
Call it $\sigma_i$.

Since $\Omega_3^{\roman{spin}}\cong0$, 
$(\widetilde{V_i}\cup D^3,$ $\sigma_i)$=
$\partial$$(W, \tau)$ 
for a compact oriented 4-manifold $W$ with a spin structure $\tau$.

The following is known. See e.g. [3]. 
We can suppose that there exists a handle decomposition 

$W$=
$\widetilde{V_i}\x[0,1]$ $\cup$ (4-dimensional 2-handles $\widetilde{h^2}$) 
$\cup$ $D^3\x[0,1]$ $\cdot\cdot\cdot\cdot\cdot(*)$. 

From now on, we consider under the condition of Claim 2. 

Take the collar neighborhood $M\x[0,1]$ of $M$ in $X$. 
Take $\widetilde{V_i}\x[0,1]$ in $M\x[0,1]$. 
Suppose the embedding is level preserving. 

We attach the 4-dimensional 2-handles  
$\widetilde{h^2}$ in $(*)$ 
to the submanifold $\widetilde{V_i}\x[0,1]$ ($\subset M\x[0,1]\subset X$) 
by the attaching map in $(*)$.
The attach parts $\widetilde{h^2}\cap(\widetilde{V_i}\x[0,1])$ 
are in $\widetilde{V_i}\x\{1\}$$\subset M\x\{1\}=M$.

By Claim 2, we can suppose that 
the attach parts of $\widetilde{h^2}$ do not intersect 
$\widetilde{V_a}$ $\cap$ $\widetilde{V_b}$.  
(Recall Claim 2.1. If necessary, 
we can move the attach parts of $\widetilde{h^2}$ by isotopy  
so that they do not intersect with $\widetilde{V_a}$ $\cap$ $\widetilde{V_b}$.)

Then we can attach the 5-dimensional 2-handles  
$\widetilde{h^2}\x [-1,1]$ to $X$ naturally. 
The tubular neighborhood of the attach part of 
 $\widetilde{h^2}$ in $M$ is the attach part of 
the 5-dimensional 2-handle  $\widetilde{h^2}$ in $M$.  

Put $\widetilde Y$=$X$ $\cup$ (the 5-dimensional 2-handles).

By the definition, 
the attaching maps of the 5-dimensional 2-handles are spin preserving diffeomorphism. 
Hence 
 $\widetilde Y$= $\overline{(\natural^q S^2\x B^3)-B^5}$ 
for a non-positive integer $q$.

Put $Y$= $\widetilde Y$ and take $f$ to be 
$f(S^2_i\x[0,1])$=$(\partial V_i)\x[0,1]$, 
where $\partial V_i\x[0,1]$ is defined in the proof of Lemma 1 ($i=a,b$).

Hence Lemma 3 holds by the construction of them.

Note that the following holds under the condition of Lemma 3.  

\proclaim{Claim}

$S^1_a$ is in the interior of $D^3_b$.  
$S^1_a$ in $D^3_b$ is the trivial knot. 

$S^1_b$ is in the interior of $D^3_a$. 
$S^1_b$ in $D^3_a$ is the trivial knot. 
\endproclaim

By using the collar neighborhood of $N$ in $Y$, 
we have the following. 
\proclaim{Claim}
There exists an immersion 
$h:$ $\widetilde{D^3_a}$ $\amalg$ $\widetilde{D^3_b}$ 
$\i Y\cong\overline{\natural^q S^2\x B^3-B^5}$  
with the following properties, 
where $\widetilde{D^3_i}$ is diffeomorphic to the 3-disc.  
\roster
\item
$h$ is transverse to $S^4_0$. 

$h(\widetilde{D^3_i})$ $\cap\partial S^4_0$=$h(\partial\widetilde{D^3_i})$ 
\item
$h(\partial\widetilde{D^3_a}\cup\partial\widetilde{D^3_b})$ 
in $S^4_0$ is the 2-link $L$. 
\item
$h(\partial\widetilde{D^3_a})$ $\cap$ $h(\partial\widetilde{D^3_b})$  
is one circle $C$.  
 
$C$ in $h(\partial\widetilde{D^3_i})$ is the trivial knot ($i=a,b$). 
\endroster
\endproclaim

$Y\cong\overline{\natural^q S^2\x B^3-B^5}$ 
is embedded in $B^5$ so that $S^4_0$=$\partial B^5$.
Therefore Theorem 4.1.1 holds.

We use `multiple' handles in the proof of Lemma 1 and 3. 
By using such handles as in a similar manner, 
we prove 
Theorem 4.1 in the case when $2m\neq2$ and 
Theorem 4.2.

\head 6. Further results\endhead
In [4], the author proved the following theorem. (Theorem1.1 of [4]). 
By  Theorem 4.1.1 in this paper and Theorem 1.1 of [4],       
the following Theorem 6.1 holds. 
In order to state Theorem 1.1of [1], we need some preliminaries. 
\f {\bf Definition }
$(L_1, L_2, X_1, X_2)$ is called a 
{\it 4-tuple of links } if the following conditions (1), (2) and (3) hold.

\f (1) $L_i=(K_{i1},...,K_{im_i})$  is an oriented ordered $m_i$-component 
1-dimensional link  $(i=1,2).$  \newline

\f (2)    $m_1=m_2.$

\f (3)     $X_i$ is an oriented 3-knot. 

\f {\bf Definition }
A 4-tuple of links $(L_1, L_2, X_1, X_2)$ 
is said to be  {\it realizable } if there exists 
a smooth transverse immersion  
$f:S^3_1\coprod S^3_2 \looparrowright S^5$  
satisfying the following conditions (1) and (2). 

\f(1)   $f\vert S^3_i$ is a smooth embedding 
and defines the  3-knot $X_i (i=1,2)$ in $S^5$.  

\f(2)   For  $C=f(S^3_1)\cap f(S^3_2)$, 
the inverse image  $f^{-1}(C)$ in $S^3_i$ defines the  1-link 
$L_i (i=1,2).$ Here, the orientation of $C$ is induced 
naturally from the preferred orientations of $S^3_1, S^3_2,$ and $S^5$,   
 and an arbitrary order is given to the components of $C$.

The following theorem characterizes the realizable 4-tuples of links.

\proclaim{Theorem 1.1 of [4]}
A 4-tuple of links $(L_1, L_2, X_1, X_2)$ is realizable  if and only if 
$(L_1, L_2, X_1, X_2)$ satisfies 
one of the following conditions i) and  ii). 

i) Both $L_1$ and $L_2$ are proper links, and 

$$\roman{Arf}(L_1) = \roman{Arf}(L_2).$$ 

ii) Neither $L_1$ nor $L_2$ 
is proper, and 
$$ lk(K_{1j}, L_1-K_{1j})
\equiv
  lk(K_{2j}, L_2-K_{2j}) \quad \roman{mod\hskip2mm  2}
 \quad \roman{for \hskip2mm all}\hskip2mm j. $$ 

\endproclaim

In order to state Theorem 6.1, we need some preliminaries. 

In Problem 3.1  in the case when $2m=2$, 
we have two 1-dimensional links, 
$D^3_1$ $\cap$ $D^3_2$ in  $D^3_1$ and $D^3_1$ $\cap$ $D^3_2$ in  $D^3_2$.

($P_1,P_2$) is called a {\it pair of 1-links} if 
$P_i$ is a 1-link $(P_{11},...,P_{1m_i})$
and $m_1$=$m_2$.

A pair of 1-links ($P_1,P_2$) is said to be {\it realizable} 
if there exist $D^3_1$ and $D^3_2$ as above such that 
the 1-link $P_i$ is equivalent to the 1-link 
$D^3_1$ $\cap$ $D^3_2$ in  $D^3_i$, 
where we give the 1-link $D^3_1$ $\cap$ $D^3_2$ in  $D^3_i$ an order.  

\proclaim{Theorem 6.1}
A pair of 1-links ($P_1,P_2$) is realizable if and only if 
one of the following conditions (1) and (2) holds. 
 
\roster
\item
Both $P_1$ and $P_2$ are proper links, and 
$$\roman{Arf}(P_1) = \roman{Arf}(P_2).$$ 
\item
 Neither $P_1$ nor $P_2$ is proper, and 
$$ lk(P_{1j}, P_1-P_{1j})
\equiv
  lk(P_{2j}, P_2-P_{2j}) \quad \roman{mod\hskip2mm  2}
 \quad \roman{for \hskip2mm all}\hskip2mm j. $$
\endroster
\endproclaim

\enddocument